\documentclass[a4paper]{article}
\usepackage[french]{babel}
\usepackage{amsfonts}
\usepackage{amsmath}
\usepackage{amsmath,amstext, amsthm, amssymb}

\font \tenmat = msbm10 \font \sevenmat = msbm7 \font \fivemat =
msbm5 \newfam \matfam \textfont \matfam = \tenmat \scriptfont
\matfam = \sevenmat \scriptscriptfont \matfam = \fivemat \def \mat
{\fam \matfam \tenmat}

\def\NN{{\mat N}} \def\ZZ{{\mat Z}}  \def\RR{{\mat R}} \def\CC{{\mat C}}

\def\Log{\mathrm{Log\,}}

\newtheorem{thm}{Th\'eor\`eme}

\newtheorem{cor}{Corollaire}

\title{La s\'erie enti\`ere $1+\frac z{\Gamma(1+i)}+\frac{z^2}{\Gamma(1+2i)}+\frac{z^3}{\Gamma(1+3i)}+...$ poss\`ede une fronti\`ere naturelle~!}
\author{Changgui ZHANG\footnote{Laboratoire P. Painlev\'e (UMR -- CNRS 8524), UFR Math.,
Universit\'e de Lille 1, Cit\'e scientifique, 59655 Villeneuve d'Ascq cedex, France. Courriel : Changgui.zhang@math.univ-lille1.fr}}
\date{}

\begin{document}

\maketitle

{\bf R\'esum\'e --}
Les s\'eries lacunaires sont les exemples les plus classiques des s\'eries enti\`eres qui ne peuvent pas se prolonger au del\`a de leurs cercles de convergence (\cite[\S~93-94, p.~372-383]{Di}, \cite[\S 7.43, p.~223]{Ti}, ...). Dans la pr\'esente Note, nous \'etudions une famille de s\'eries enti\`eres, non lacunaires, ayant pour coefficients des valeurs prises par la fonction Gamma sur des lignes verticales. Nous expliquons comment les repr\'esenter en termes de s\'eries de Dirichlet lacunaires, ce qui nous permet de conclure \`a l'existence de leur fronti\`ere naturelle. Li\'es au comportement \og al\'eatoire\fg\ de la fonction Gamma sur toute ligne verticale, les r\'esultats ainsi obtenus verront \'egalement des explications dans notre travail en cours sur des \'equations aux $q$-diff\'erences-diff\'erentielles, dites \og type pantagraphe \fg \ (voir \cite{KM} pour l'instant).

\bigskip

{\bf Abstract -- The Power Series $1+\frac z{\Gamma(1+i)}+\frac{z^2}{\Gamma(1+2i)}+\frac{z^3}{\Gamma(1+3i)}+...$ Has a Natural Boundary --}
The lacunary series are the most classical examples among all the power series whose circle of convergence constitutes a natural boundary (\cite[\S~93-94, p.~372-383]{Di}, \cite[\S 7.43, p.~223]{Ti}, ...). In this Note, we study a family of non-lacunary power series whose coefficients are given by means of values of the Gamma function over vertical line. We explain how to transform these series into lacunary Dirichlet series, which allows us to conclude the existence of their natural boundary. Our results, which illustrate in what manner the Gamma function may have a  unpredictable behaviour on any vertical line, may also be partially understood in the framwork of our forthcoming work on a class of differential $q$-difference equations, namely, on pantagraph type equations (see \cite{KM} for instance). 

\bigskip

\bigskip

\centerline{------}

\bigskip

\bigskip

Dans un travail sur des \'equations lin\'eaires fonctionnelles aux $q$-diff\'erences et diff\'erentielles (voir, pour l'instant, \cite{KM}), nous rencontrons une famille de s\'eries de Laurent qui \'evoquent des valeurs de la fonction Gamma sur une ligne verticale. Le but de la pr\'esente Note est d'\'etudier leur prolongement analytique et de faire remarquer l'existence de la fronti\`ere naturelle.

Bien que la question de la coupure analytique, en tant que sujet de recherche g\'en\'eral, semble devenue \og d\'emod\'ee\fg, nos s\'eries, \`a coefficients \og explicites\fg,  fournissent n\'eanmoins des exemples  \og concrets et naturels\fg\ de s\'eries non lacunaires qui ne peuvent se prolonger au del\`a du disque de convergence. Selon une id\'ee g\'en\'eralement re\c cue, depuis le travail d'E. Fabry \cite{Fa}, chez une s\'erie \`a fronti\`ere naturelle, les coefficients de cette derni\`ere para\^\i tront comme \'etant \og arbitrairement donn\'es\fg.  Ainsi nos s\'eries permettent-elles d'illustrer \`a quel point est \og arbitraire\fg\ le comportement de $\Gamma(a+ib)$ lorsque $b$ tend vers l'infini suivant une progression arithm\'etique. Voir l'article \cite{AR} pour une \'etude sur la convexit\'e de $\log \Gamma(z)$ dans le plan complexe et, en particulier, sur une droite verticale de celui-ci.

Le reste de l'article comprendra deux paragraphes : le th\'eor\`eme \ref{thm:Psiuvz+}, notre principal r\'esultat, sera \'enonc\'e dans le premier et sa d\'emonstration sera faite dans le dernier.

\section{Notations et \'enonc\'es}\label{section:1}

Les notations $\NN$, $\ZZ$, $\RR$, $\CC$ sont standard.

\medskip

{\bf (1.1)}  Soit $(u,v)\in\CC\times\RR$ tel que $u\notin -\NN+\frac{2vi}{\pi}\ZZ$ et consid\'erons la s\'erie de Laurent not\'ee $\Psi(u,v,z)$, de la variable $z$, d\'efinie par
$$
\Psi(u,v,z)=\sum_{n\in\ZZ}\Gamma(u+\frac{2ivn}{\pi})\,z^n\,.
$$
D'apr\`es la formule de Stirling, on a ({\it cf} \cite[Corollary 1.4.4, p.~21]{AAR}):
$$
\Gamma(u+\frac{2ivn}{\pi})=O(n^{u-\frac{1}{2}}\,e^{-\vert vn\vert})
$$
lorsque l'indice $n$ tend vers $\pm \infty$; on en d\'eduit que, si $v\not=0$, la s\'erie $\Psi(u,v,z)$ d\'efinit une fonction analytique dans la couronne ${\cal C}_v$ avec ceci:
$${\cal C}_v:=\{z\in\CC: e^{-\vert v\vert}<\vert z\vert<e^{\vert v\vert}\}.$$

Soit $\CC^+$ le demi-plan $\Re u>0$ de $\CC$ et $\RR^+$ la demi-droite $v>0$ de $\RR$.

\begin{thm}\label{thm:Psiuvz+}
Pour tout couple $(u,v)\in\CC^+\times\RR^+$, la s\'erie $\Psi(u,v,z)$ admet $\partial {\cal C}_v$ pour fronti\`ere naturelle.
\end{thm}

\medskip

{\bf (1.2)} Ce th\'eor\`eme sera d\'emontr\'e dans le parapraphe \ref{section:2}; avant de faire ceci, nous nous contenterions de donner quelques cons\'equences du r\'esultat. 

D'abord, notons que la fonction $\Psi(u,v,z)$ satisfait aux relations suivantes:
$$
\Psi(u+1,v,z)=(\frac{2vi}\pi\partial_z+u)\Psi(u,v,z)\,,\quad
\Psi(u,v,\frac1z)=\Psi(u,-v,z)\,,
$$
lesquelles nous permettent d'\'etendre le th\'eor\`eme \ref{thm:Psiuvz+} de la mani\`ere suivante.

\begin{cor}\label{cor:Psiuvz}
Pour tout couple $(u,v)\in\CC\times \RR^*$, si $u\notin-\NN+\frac{2vi}{\pi}\ZZ$, alors la s\'erie de Laurent $\Psi(u,v,z)$ admet $\partial {\cal C}_v$ pour fronti\`ere naturelle.
\hfill $\square$
\end{cor}

En consid\'erant $\Psi(u,v,z)$ comme \'etant la somme d'une s\'erie enti\`ere de la variable $z$ avec une s\'erie enti\`ere de $\frac 1z$, nous d\'eduisons ais\'ement du corollaire \ref{cor:Psiuvz} le

\begin{cor}\label{cor:uvz}
Soit $(u,v)\in\CC\times \RR^*$ tel que $-u\notin\NN+\frac{2iv}{\pi}\NN$. La s\'erie enti\`ere $\sum_{n\ge 0}\Gamma(u+\frac{2ivn}{\pi})z^n$ admet le cercle $\vert z\vert =e^{\vert v\vert}$ pour fronti\`ere naturelle.
\hfill $\square$
\end{cor}

\medskip

{\bf (1.3)} D'apr\`es la formule du compl\'ement de $\Gamma$, on a la relation
$$
\Gamma(u+\frac{2ivn}{\pi})\Gamma(1-u-\frac{2ivn}{\pi})=\frac{2\pi i}{e^{\pi ui-{2vn}}-e^{-\pi ui+{2vn}}}\,,
$$
ou encore:
$$
\frac{2\pi i}{\Gamma(u+\frac{2ivn}{\pi})}=\left({e^{\pi ui-{2vn}}-e^{-\pi ui+{2vn}}}\right)\,\Gamma(1-u-\frac{2ivn}{\pi})\,,
$$
laquelle, combin\'ee avec le corollaire \ref{cor:uvz}, nous conduit au r\'esultat suivant, jusitifiant ainsi le titre de l'article.

\begin{cor}\label{cor:uv1}
Soit $(u,v)\in\CC\times \RR^*$. La s\'erie enti\`ere $\sum_{n\ge 0}\frac{z^n}{\Gamma(u+\frac{2ivn}{\pi})}$ admet le cercle $\vert z\vert =e^{-\vert v\vert}$ pour fronti\`ere naturelle.
\hfill $\square$
\end{cor}

\medskip

{\bf (1.4)} Les r\'esultats pr\'esent\'es ci-dessus trouveront une explication partielle dans un ordre d'id\'ees assez proche de celui du travail d'E. Fabry \cite{Fa}. En effet, le couple $( u,v)$ \'etant fix\'e dans $\CC^+\times\RR^+$, si l'on note $\phi(n)=\arg\Gamma(u+\frac{2ivn}{\pi})$ $\bmod\ 2\pi$, la formule de Stirling implique que, lorsque $n$ augmente ind\'efiniment, $\phi(n)$ est d'ordre $n\ln n$ et
$$
\phi(n+1)-\phi(n)=\frac{2v}{\pi}\,\ln\bigl(\frac{2vn}{\pi}\bigr)+O\bigl(\frac{1}{n}\bigr)\quad\bmod\ 2\pi.
$$

\medskip

{\bf (1.5)}
Mentionnons enfin qu'il serait plausible d'\'etendre les r\'esultats de (1.1)-(1.3) aux s\'eries pour lesquelles les coefficients $\alpha_n$ sont donn\'es par l'une des m\'ethodes suivantes.

\medskip

(A)  $\alpha_n=R(\Gamma(u+\frac{2ivn}{\pi}))$, $R$ \'etant une fraction rationnelle d'une variable.

\medskip

(B) $\alpha_n=(a_{1,n}...a_{\ell,n})/(b_{1,n}...b_{m,n})$, o\`u $\ell$, $m\in\NN$, $\ell+m>0$ et o\`u chaque facteur $a_{j,n}$ ou $b_{j,n}$ est de la forme $\Gamma(u_j+\frac{2iv_jn}{\pi})$, les $(\ell+m)$ couples $(u_j, v_j)$ \'etant suppos\'es \og g\'en\'eriquement non li\'es\fg.

\medskip

(C) Combiner les m\'ethodes (A) et (B).

\section{D\'emonstration du th\'eor\`eme \ref{thm:Psiuvz+}}\label{section:2}

Elle se fait au moyen d'une repr\'esentation de $\Psi(u,v,z)$ en termes de s\'erie de Dirichlet, obtenue depuis l'int\'egrale eul\'erienne de $\Gamma$ et du th\'eor\`eme des r\'esidus.

Dans ce paragraphe, nous rappelons que $\Re u>0$ et $v>0$.

\medskip

{\bf (2.1)} Soit $L_+$ et $L_-$ deux demi-droites partant de l'origine dans le premier quadrant et dans le quatri\`eme quadrant respectivement, c'est-\`a-dire, $L_+=]0,\infty e^{i\epsilon}[$ et $L_-=]0,\infty e^{-i\epsilon'}[$ avec $\epsilon$, $\epsilon'\in]0,\frac{\pi}{2}[$; afin d'alleger l'exposition, nous supposerons $\epsilon'=\epsilon$. A l'aide de l'int\'egrale eul\'erienne de $\Gamma$, on \'ecrit:
$$
\Psi(u,v,z)=\sum_{n\ge 0}\int_{L_+}e^{-t}\,t^u\,\bigl(t^{\frac{2iv}{\pi}}z\bigr)^n\frac{dt}{t}+\sum_{n<
0}\int_{L_-}e^{-t}\,t^u\,\bigl(t^{\frac{2iv}{\pi}}z\bigr)^n\frac{dt}{t}\,,
$$
o\`u $t^\alpha=e^{\alpha\Log t}$, $\Log$ d\'esignant la branche principale du logarithme dans le plan complexe ou, plus exactement, dans la surface de Riemann associ\'ee.

Supposons que $z\in\CC$ satisfasse aux conditions suivantes :
$$
\sup_{t\in L_+}\bigl\vert zt^{\frac{2iv}{\pi}}\bigr\vert<1,\qquad
\sup_{t\in L_-}\bigl\vert zt^{\frac{2iv}{\pi}}\bigr\vert>1\,,
$$
ou, de mani\`ere \'equivalente,
$$
e^{-\frac{2v}{\pi}\epsilon}<\vert z\vert<e^{\frac{2v}{\pi}\epsilon},\qquad
\epsilon\in]0,\frac{\pi}{2}[\,;
$$
sous ses conditions, et par un argument standard de type Fibini, on arrive \`a la repr\'esentation int\'egrale suivante:
$$
\Psi(u,v,z)=\int_{L_+-L-}\frac{e^{-t}\,t^u}{1-zt^{\frac{2iv}{\pi}}}\,\frac{dt}{t}\,.\leqno(I)
$$

\medskip

{\bf (2.2)} Soit $z\in{\cal C}_v$ et choisissons $\epsilon$ suffisammenet proche de $\frac{\pi}{2}$ dans $]0,\frac{\pi}{2}[$ de fa\c con \`a avoir l'int\'egrale (I) pour $\Psi(u,v,z)$; supponsons, en plus, que $\arg z\in]-{\pi},\pi]$. En annulant le d\'enomi\-nateur $1-zt^{\frac{2iv}{\pi}}$ sous le signe int\'egral de (I), on trouve que les p\^oles, simples, de l'int\'egrand dans le demi-plan \`a droite constituent la \og spirale\fg $t_z{\tilde v}^\ZZ$, avec
$$
t_z=z^{-\frac{\pi i}{2v}}=e^{-\frac{\pi i}{2v}\Log z},\qquad
\tilde v=e^{\frac{\pi^2}v}\in]1,+\infty[\,;
$$
on remarquera que cette spirale se situe entre $L_+$ et $L_-$.

Ceci \'etant, le th\'eor\`eme des r\'esidus appliqu\'e \`a l'int\'egrale (I) nous conduit \`a l'expression suivante:
$$
\Psi(u,v,z)=\frac{\pi^2e^{-\frac{u\pi i}{2v}\Log z}}{v}\sum_{k\in\ZZ}{\tilde v}^{ku}\,e^{-t_z{\tilde v}^k}\,,\leqno(D)
$$
dans laquelle la s\'erie du second membre est une \og s\'erie de Dirichlet\fg; voir \cite[Chapter IX, \S 8, p. 432-440]{SZ}.

\medskip

{\bf (2.3)} Soit $u\in\CC^+$, $q>1$ et consid\'erons les s\'eries suivantes:
$$
S_-(u,q,\zeta)=\sum_{k<0}q^{ku}e^{-\zeta q^k}\,,\qquad
S_+(u,q,\zeta)=\sum_{k\ge 0}q^{ku}e^{-\zeta q^k}\,.\leqno(S)
$$
Il est clair que $S_-(u,q,\zeta)$ d\'efinit une fonction enti\`ere de la variable $\zeta$.
D'autre part,  $S_+(u,q,\zeta)$ d\'efinit une fonction analytique dans le demi-plan $\Re\zeta>0$ et l'axe imaginaire $\Re\zeta=0$ en forme une fronti\`ere naturelle: voir l'appendice \`a la fin de l'article. On terminera ainsi la preuve du th\'eor\`eme \ref{thm:Psiuvz+} gr\^ace \`a l'expression (D), avec $q=\tilde v$ et $\zeta=t_z$.\hfill $\square$

\medskip

{\bf (2.4)} Observons enfin que, dans (S), on a les relations fonctionnelles 
$$q^uS_-(u,q,q\zeta)-S_-(u,q,\zeta)=e^{-\zeta}\,,\qquad
q^uS_+(u,q,q\zeta)-S_+(u,q,\zeta)=-e^{-\zeta}\,;
$$
ou encore, en posant $S(u,q,\zeta)=S_-(u,q,\zeta)+S_+(u,q,\zeta)$:
$$
q^uS(u,q,q\zeta)=S(u,q,\zeta).
$$
Ces relations repr\'esentent des \'equations lin\'eaires aux $q$-diff\'erences singuli\`eres r\'eguli\`eres (ou dites fuchsiennes) en $\zeta=0$; voir \cite{DRSZ} et \cite{Zh}.

\medskip

{\bf Appendice  --} Par analogie avec la th\'eorie des s\'eries enti\`eres lacunaires, voici un \'enonc\'e concernant la s\'erie de Dirichlet. {\it Soit $\lambda_1<\lambda_2<...$ une suite de r\'eels strictement croissante telle que $\lambda_{n+1}-\lambda_n \to \infty$ pour $n\to\infty$. Soit $A(s)=\sum_{n\ge 1}a_ne^{-\lambda_n s}$ une s\'erie de Dirichlet, avec $a_n\in\CC$; on notera $\sigma=\sigma_c(A)$ son abscisse de convergence. Si $\sigma$ est finie, alors l'axe $\Re s=\sigma$ constitue une coupure analytique pour la s\'erie $A$.} 
\medskip

{\it Remerciement --} L'Auteur  tient \`a remercier  Professeur Herv\'e Queff\'elec de lui avoir \og prouv\'e \`a la main\fg\ l'\'enonc\'e inclus dans l'appendice et remercier Professeur Jean-Fran\c cois Burnol pour ses commentaires et remarques.

\end{document}